\documentclass[11pt]{article}

\usepackage{setspace}
\usepackage{amssymb,amsmath,amsthm}
\usepackage{amsfonts}
\usepackage{calc}
\usepackage{verbatim}
\usepackage{epsfig}
\usepackage{graphicx}
\usepackage{graphics}
\usepackage{xcolor}\usepackage[hmargin=3cm,vmargin=3.5cm]{geometry}
\usepackage{chngcntr}
\counterwithin{figure}{section}
\newtheorem{defn}{Definition}[section]

\newtheorem{prop}[defn]{Proposition}
\newtheorem{theorem}[defn]{Theorem}
\newtheorem{cor}[defn]{Corollary}

\newtheorem{ques}[defn]{Question}

\newcommand{\CC}{\mathbb{C}}
\newcommand{\ee}{\mathrm{e}}
\newcommand{\ii}{\mathrm{i}}
\newcommand{\drop}{\smallsetminus}
\newcommand{\RR}{\mathbb{R}}
\newcommand{\showon}{\begin{eqnarray*}}
\newcommand{\showoff}{\end{eqnarray*}}
\newcommand{\dd}{\mathrm{d}}
\newcommand{\Disc}{\mathrm{Disc}}

                                {\null\hfill$\Box$\par\medskip\medskip\medskip\medskip}

\makeatletter
\renewcommand{\pod}[1]{\mathchoice
  {\allowbreak \if@display \mkern 18mu\else \mkern 8mu\fi (#1)}
  {\allowbreak \if@display \mkern 18mu\else \mkern 8mu\fi (#1)}
  {\mkern4mu(#1)}
  {\mkern4mu(#1)}
}

\title{On the imaginary parts of chromatic roots}

\author{Jason I. Brown\\
Department of Mathematics and Statistics\\
Dalhousie University\\
Halifax, Nova Scotia, Canada B3H 3J5\\
jason.brown@dal.ca
\and
David G. Wagner\\
Department of Combinatorics and Optimization, University of Waterloo,\\ Waterloo, Canada N2L 3G1\\ dgwagner@uwaterloo.ca}

\begin{document}

\maketitle

%\biboptions{sort&compress}

\begin{abstract}
While much attention has been directed to the maximum modulus and maximum real part of chromatic roots of graphs of order $n$ (that is, with $n$ vertices), relatively little is known about the maximum {\em imaginary part} of such graphs. We prove that the maximum imaginary part can grow linearly in the order of the graph. We also show that for any fixed $p \in (0,1)$, almost every random graph $G$ in the Erd\"{o}s-R\'{e}nyi model has a non-real root.   
\end{abstract}

\section{Introduction}
A (vertex) $k$-colouring of a (finite, undirected, simple) graph $G$ is a function $f:V(G) \rightarrow \{1,\ldots,k\}$ such that no two adjacent vertices receive the same colour, that is, if $uv$ is an edge of $G$, then $f(u) \neq f(v)$. The function $\pi(G,k)$ that for all nonnegative integers $k$ counts the number of $k$-colourings of $G$, is well known to be a polynomial function of $k$, and its extension to all complex numbers $x$ is called the {\em chromatic polynomial} of $G$. There is likely no better studied graph polynomial than the chromatic polynomial, with interest initiated by Birkhoff in his work on the famous Four Colour Conjecture -- whether every planar graph can be coloured with four colours. The research literature on the topic is vast -- see \cite{read}, \cite{readtutte} and \cite{dongbook} for a recent survey.

The Four Colour Theorem \cite{fct1,fct2} is equivalent to stating that $4$ is never a root of the chromatic polynomial of a planar graph, and the nature and location of roots of chromatic polynomials ({\em chromatic roots}) has been of great interest. There are no negative real roots (as the coefficients of a chromatic polynomial alternate in sign).
% while an easy argument shows that the real roots of a graph of order $n$ (that is, with $n$ vertices) is bounded above by $n-1$: if $c_{k}$ denotes the number of ways to partition the vertices of a graph $G$ of order $n$ into $k$ {\em independent sets} (i.e., subsets of vertices that contain no edges), and $(x)_{\downarrow k} = x(x-1)(x-2) \cdots (x-k+1)$ denote the {\em $k$--th falling factorial} of $x$, then 
%\[ \pi(G,x) = \sum_{k} c_{k}(x)_{\downarrow k}.\]
%Each $c_{k}$ is nonnegative with $c_{n}$ certainly positive, and $(x)_{\downarrow k}%$ is always positive for $x > k$, completing the argument. 
While it is known \cite{jackson,thomassen} that the closure of the set of real chromatic roots is the set $\{0,1\} \cup 
[32/27,\infty)$, the closure of the set of complex chromatic roots is in fact the whole complex plane \cite{sokaldense}.

Various results are known about the maximum modulus of chromatic roots of a graph, with regard to the order and size (that is, the number of edges). For example, the chromatic roots of graphs of order $n$ and size $m$ are known to be in the disks $|z| < 8 \Delta < 8n$ \cite{sokalbound} and $|z-1| \leq m-n+1$  \cite{brownmonomial}, where $\Delta$ is the maximum degree of a vertex of $G$. On the other hand, there are chromatic roots of modulus at least $\displaystyle{\frac{m-1}{n-2}}$ \cite{brownbound}. As every complete graph of order $n$ has a chromatic root at $n-1$, the rate of growth of 
\[ \mathrm{maxmod}(n) = \mathrm{max} \{ |z| : z \mbox{ is a chromatic root of a graph of order } n \}\]
is linear, that is, there are positive constants $C_{1}$ and $C_{2}$ such that $C_{1}n \leq \mathrm{maxmod}(n) \leq C_{2}n$. The same result holds for the maximum {\em real part} of a chromatic root of a graph of order $n$,
\[ \mathrm{maxreal}(n) = \mathrm{max} \{ \Re(z) : z \mbox{ is a chromatic root of a graph of order } n \}\]
-- the function also grows linearly, for the same reasons.

However, what can be said about the growth rate of 
\[ \mathrm{maximaginary}(n) = \mathrm{max} \{ \Im(z) : z \mbox{ is a chromatic root of a graph of order } n \},\]
the maximum {\em imaginary part} of a chromatic root of order $n$? Very little is known. In \cite{brownbound} it was shown that the maximum imaginary part of a chromatic root of the complete bipartite graph $K_{\lfloor n/2 \rfloor,\lceil n/2 \rceil}$ of order $n$ has a chromatic root with imaginary part $\Omega(\sqrt{n})$. Of course, from the maximum modulus of a chromatic root of a graph of order $n$ being at most $8 \Delta$, the growth rate of the maximum imaginary part of a chromatic root is no more than linear. Alan Sokal (private communication) has suggested that, via computations, the maximum imaginary part of the complete bipartite graph $K_{n,n}$ seems to be about $.7239685$ times the order of the graph. However, a rigorous argument is elusive. What is the true rate of growth of the maximum imaginary part of a chromatic root of a graph of order $n$?

Another question concerns chromatic roots for random graphs. Only some computational results are to be found \cite{bussel}. Many graphs, including forests and chordal graphs (including complete graphs) have only real chromatic roots, while others (such as complete bipartite graphs) do not. We ask: do almost all graphs (in the Erd\"{o}s-R\'{e}nyi model, for fixed edge probability $p \in (0,1)$) have all real roots? We show that for any fixed $p \in (0,1)$, almost all graphs, in fact, have a non-real chromatic root. 

%Finally we turn to a related graph polynomial, the {\em $\sigma$-polynomials}. For a graph $G$, the $\sigma$-polynomial of $G$ is 
%\[ \sigma(G,x) = \sum a_{i}x^{i},\] 
%where $a_{i} = a_{i}(G)$ is the number of partitions of the vertices of $G$ into $n-i$ independent sets (so $a_{0} = n$ and $a_{1} = {{n} \choose {2}} - m$, where $m$ is the number of edges in $G$).
%Almost all graphs of order at most $9$ have all real $\sigma$-roots (a $\sigma$-root is a root of a $\sigma$-polynomial of some graph). There are only two graphs of order at most $8$ that have a non-real $\sigma$-root, and only $42$ such graphs of order $9$ \cite{brentiroylewagner}. As the $\sigma$-polynomial is multiplicative over joins of graphs (the {\em join} of two vertex disjoint graphs is the graph formed from their disjoint union by adding all edges between the graphs), we can trivially build larger graphs with the same non-real $\sigma$-roots. However, there are no nontrivial infinitely families of graphs (not built from smaller ones via the join operation) with non-real $\sigma$-roots. We shall present here some new infinite families of graphs with non-real $\sigma$-roots. 

Our approach to both problems involves the well known \textit{Gauss-Lucas Theorem}:

\begin{theorem}[Gauss-Lucas] \label{GL}
Let $f(z)\in\CC[z]$ be a nonconstant polynomial with complex coefficients, and
let $f'(z) = \dd f(z)/\dd z$ be the derivative of $f(z)$.  Then all
the roots of $f'(z)$ lie in the convex hull of the set of roots of $f(z)$ in $\CC$.
\end{theorem}

A simple but important consequence is the following:

\begin{cor} If some nonzero iterated derivative of a polynomial $f(x)$ with complex coefficients has a root with imaginary part $b > 0$, then $f(x)$ has a non-real root as well, with imaginary part at least $b$.
\end{cor}

A proof of Theorem \ref{GL} can be found in \cite{prasolov}.
Informally, the proof is very simple.  Let $f(z)\in\CC[z]$ be a polynomial of degree $d\geq 1$,
and let $\xi_1,...,\xi_d$ be the roots of $f(z)$ (which are not necessarily distinct).
Then $f'(z)/f(z) = \sum_{j=1}^d (z-\xi_j)^{-1}$ as rational functions in $\CC(z)$.
Let $K$ be the convex hull of
$\{\xi_1,...,\xi_d\}$, and consider any $w\in\CC\drop K$.  Note that $f(w)\neq 0$.  There is a line $L$ in the
complex plane with $w$ on one side of $L$ and $K$ on the other side of $L$.  Let $\theta\in(-\pi,\pi]$ be either
of the angles such that the line $\ee^{\ii\theta}\RR$ is perpendicular to $L$.  Then all of the real numbers
$\Re(\ee^{-\ii\theta}(w-\xi_j))$ for $1\leq j\leq d$ are nonzero and have the same sign.  It follows
that $f'(w)/f(w)\neq 0$, so that $f'(w)\neq 0$.

We use the Gauss-Lucas Theorem \ref{GL} to investigate non-real roots of chromatic polynomials of graphs.
The idea is to differentiate a polynomial $f(z)\in\CC[z]$  repeatedly until only a polynomial $g(z)$
of degree at most four remains.  Then discriminant conditions determine whether $g(z)$ has a non-real root.
By the Gauss-Lucas Theorem \ref{GL}, if $g(z)$ has a non-real root then so does $f(z)$.

When $g(z)$ is quadratic we can solve for its roots easily, and obtain a lower bound for the largest
imaginary part of a root of $f(z)$.  For quartic polynomials $g(z)$ we use the following criterion.

\begin{prop}\cite{rees} \label{disc4}
Let $g(z) = az^4 + b z^3 + c z^2 + d z + f$ be a quartic polynomial in $\CC[z]$.
Then $g(z)$ has a non-real root if
\showon
\Disc(g(z))
&=&
256a^3f^3-192a^2bdf^2-128a^2c^2f^2+144a^2cd^2f-27a^2d^4+144ab^2cf^2-6ab^2d^2f
\\
 & & -80abc^2df+18abcd^3+16ac^4f-4ac^3d^2-27b^4f^2+18b^3cdf-4b^3d^3-4b^2c^3f+b^2c^2d^2\\
&<& 0.
\showoff
\end{prop}

The applicability of this strategy depends on the fact that the first few coefficients (of highest degree)
of chromatic polynomials have relatively straightforward combinatorial meaning.  Similarly, this reasoning
can be applied to any class of polynomials for which some of the highest order terms can be determined.
The somewhat surprising fact is that, at least in the case of chromatic polynomials, this rather weak information
seems to work quite well.

\section{Linear growth of the maximum imaginary part of a chromatic root}

While the moduli and real parts of chromatic roots grow linearly in the order of a graph, the same was not known for imaginary parts. In this section we shall prove that indeed this is the case.

\begin{figure}[htp]
\begin{center}
\includegraphics[scale=0.75]{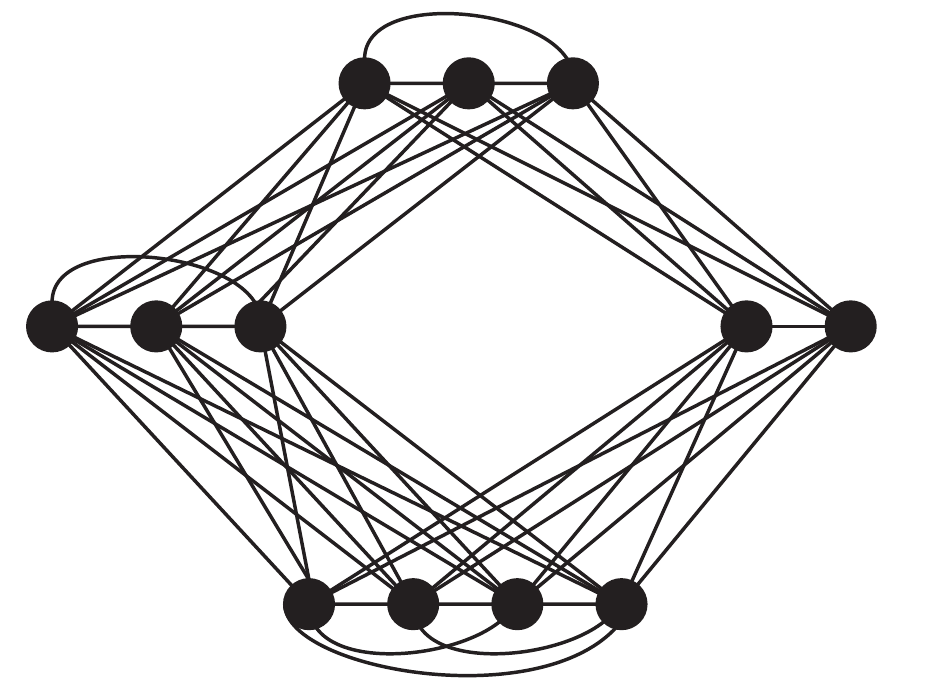}
\caption{The graph $C_{4}(3,2,4,3)$.}
\label{cycle4cliques}
\end{center}
\end{figure}

Given positive integers $a,b,c,d$, let $C_{4}(a,b,c,d)$ be the graph formed from a cycle of length $4$ by replacing the vertices in cyclic order by cliques of order $a, b, c$ and $d$ respectively (all edges are present between vertices of a clique and the two `adjacent' cliques in cyclic order). Figure~\ref{cycle4cliques} shows one such graph. 

\begin{figure}[htp]
\begin{center}
\includegraphics[scale=0.35]{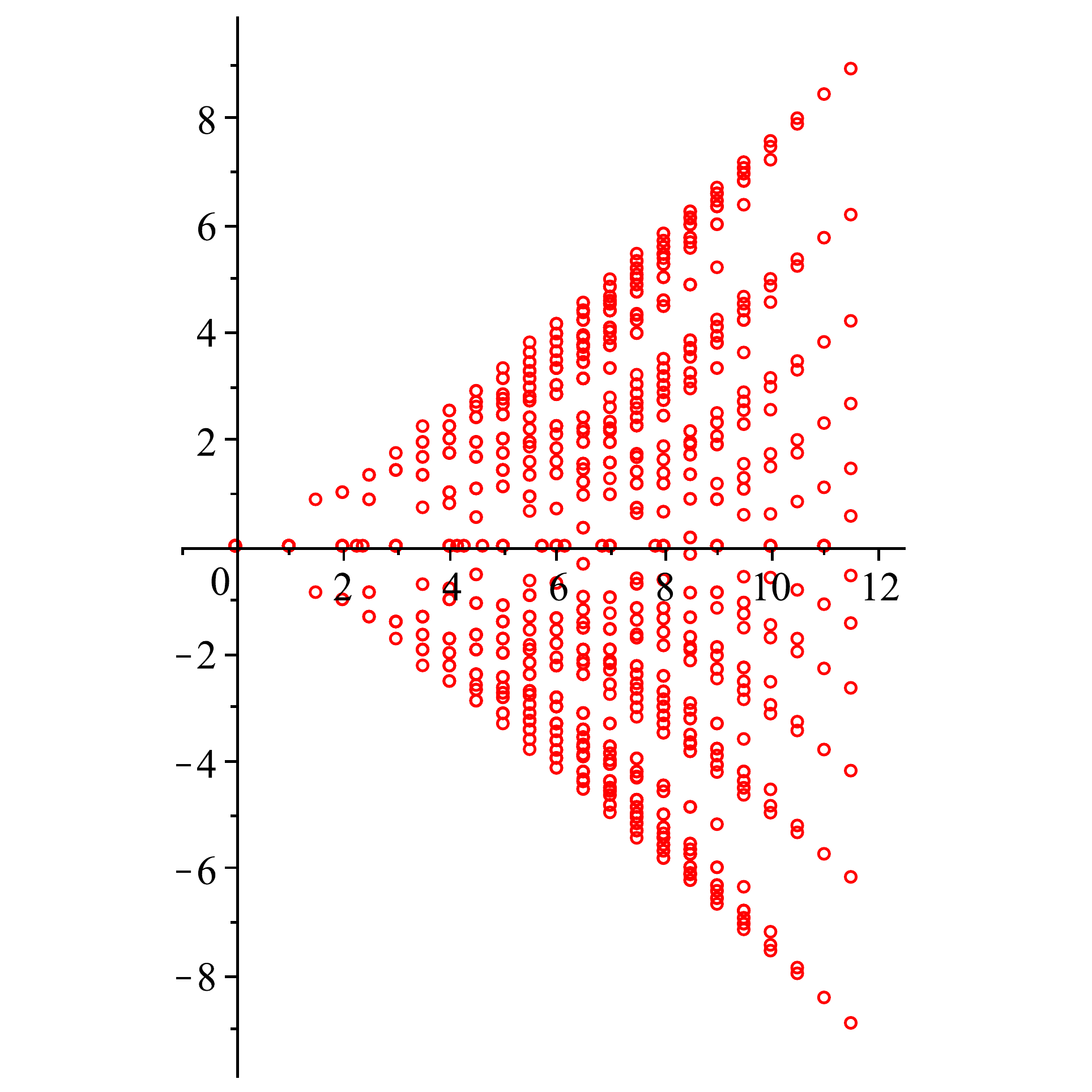}
\caption{The chromatic roots of $C_{4}(a,b,c,d)$ for $1 \leq a,b,c,d \leq 6$. Note how the roots line up in vertical stacks.}
\label{c4cliqueexample2}
\end{center}
\end{figure}

The chromatic polynomials and roots of such graphs were investigated in \cite{ringfourcliques}, where it was shown that, surprisingly, the non-real roots all have real part equal to $(a+b+c+d-1)/2$ (so that these roots line up vertically in the complex plane when the order $n = a+b+c+d$ is fixed -- see Figure~\ref{c4cliqueexample2}). 

The sharing of the same real part for the nonreal chromatic roots is interesting enough, but we are interested in the imaginary parts. To do so, we first need to travel through a sequence of polynomials related to the chromatic polynomial of $C_{4}(a,b,c,d)$ in order to discuss its chromatic roots. As in \cite{ringfourcliques}, it was observed that one could write 
\[ \pi(C_{4}(a,b,c,d),x) = \frac{(x)_{b+c}(x)_{c+d}}{(x)_{a+c}}Q_{a,b,c,d}(x),\]
where $(x)_{k} = x(x-1) \cdots (x-k+1)$ is the {\em $k$th falling factorial of $x$}, and $Q_{a,b,c,d}(x)$ is a polynomial in $x$. Moreover, if we set 
\[ p = (b+c-a-d+1)/2,~q = (c+d-a-b+1)/2, \mbox{ and } k = (b+d-a-c+1)/2,\] 
then we can express
\[ Q_{a,b,c,d}(z + (n-1)/2) = F_{a,p,q,k}(z)\]
for another polynomial $F = F_{a,p,q,k}$. Moreover, it turns out that $F$ is an even polynomial, that can be expressed as 
\[ F_{a,p,q,k}(z) = W_{a,p,q,k}(z^{2})\]
for a polynomial $W_{a,p,q,k}(z)$; in fact, $W_{a,p,q,k}(z)$ satisfies 
\begin{eqnarray}
W_{0,p,q,k}(z) & = & 1,\label{w0}\\
W_{1,p,q,k}(z) & = & z+ pq + pk + qk\label{w1}
\end{eqnarray}
and for $a \geq 2$,
\begin{eqnarray}
W_{a,p,q,k}(z) & = & (z+(a-1)(2p+2q+2k+2a-3)+pq+pk+qk)W_{a-1,p,q,k}(z) \nonumber \\  
 & & -(a-1)(p+q+a-2)(q+k+a-2)(p+k+a-2)W_{a-2,p,q,k}(z).\label{wa}
\end{eqnarray}
The roots of $W$ were then shown to be real and nonpositive, so that for every negative root $r$ of $W$, both $-\sqrt{-r} i$ and $\sqrt{-r} i$ are roots of $F$, and hence $Q$, and thus $\pi(C_{4}(a,b,c,d)$, has roots at $(n-1)/2 \pm \sqrt{-r} i$. 

With all of this out of the way, our plan is to show that when $a = b = c = d$, we can find a root $r$ of $W$ so that $-r = \Omega(n^{2})$; this will imply that the graphs $C_{4}(a,a,a,a)$ have a chromatic root with imaginary part at least $Cn$ for some positive constant $C$. We can then extend the result to all $n$ by noting that if $n = 4a+l$, for some $1 \leq l \leq 3$, then by noting that that disjoint union of a graph with isolated vertices does not change the set of chromatic roots , we find that some graph of order $n$ has imaginary part at least $C(n-3)$, and hence imaginary part at least $C^{\prime}n$ for a slightly smaller constant $C^{\prime}$ (and sufficiently large $n$).

So the question is, how large in absolute value are the roots of $W$ guaranteed to be? When $a = b= c = d$, we find from the formulas that $p = q = k = 1/2$, and that
\begin{eqnarray}
W_{0} = W_{0,1/2,1/2,1/2}(z) & = & 1,\label{w0prime}\\
W_{1} = W_{1,1/2,1/2,1/2}(z) & = & z+ 3/4\label{w1prime}
\end{eqnarray}
and for $a \geq 2$,
\begin{eqnarray}
W_{a} = W_{a,1/2,1/2,1/2}(z) & = & (z+2(a-1)a+3/4)W_{a-1}(z) - (a-1)^{4}W_{a-2}(z).\label{waprime}
\end{eqnarray}
Moreover, from this recursion, we can calculate that 
\begin{eqnarray}
W_{a} & = & z^{a} + \left( \frac{2}{3} a^{3} + \frac{1}{12} a \right) z^{a-1} + \nonumber \\
& & \left( \frac{2}{9}a^{6}-\frac{3}{5}a^{5}+\frac{5}{9}a^{4}-\frac{1}{6}a^{3}+\frac{1}{288}a^{2}-\frac{7}{480}a \right) z^{a-2} + \cdots . \label{wexpansion}
\end{eqnarray}

\begin{figure}[htp]
\begin{center}
\includegraphics[scale=0.3]{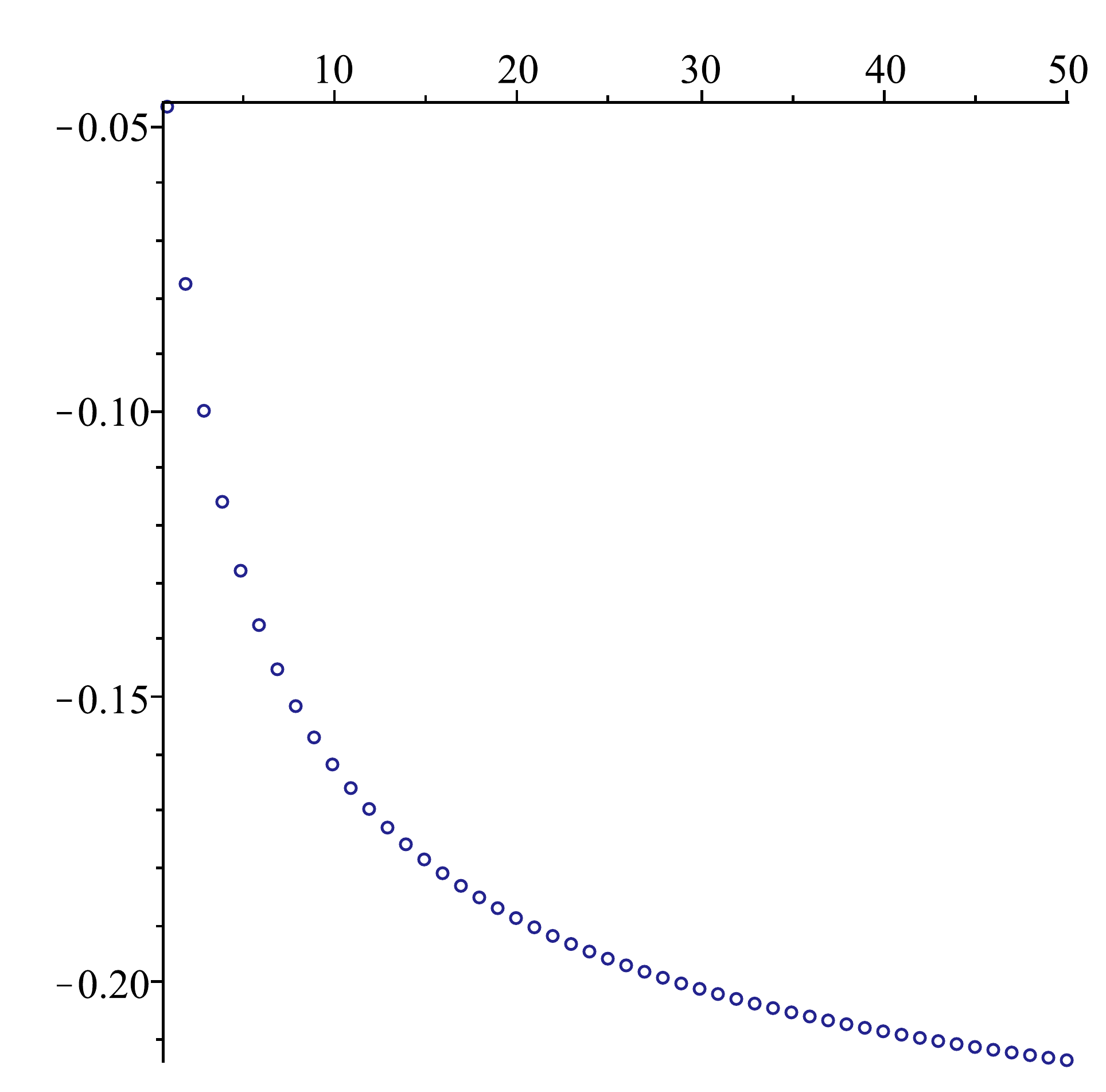}
\caption{The leftmost roots of $W_{a}$ for $a \leq 40$, divided by $n^2 = 16 a^{2}$.}
\label{warootsovernsquared}
\end{center}
\end{figure}

We are interested in the leftmost (that is, the `most negative`) root of $W_{a}$. Figure~\ref{warootsovernsquared} plots the leftmost root, divided by $n^2$, and here we see what suggests limiting behaviour. 

To get a bound on the roots of $W_{a}$, we differentiate it down $a-2$ times, until we reach a quadratic:

\[ W_{a}^{(a-2)} = \frac{a!}{2}z^{2} + (a-1)!\left( \frac{2}{3} a^{3} + \frac{1}{12} a \right) z + (a-2)!\left( \frac{2}{9}a^{6}-\frac{3}{5}a^{5}+\frac{5}{9}a^{4}-\frac{1}{6}a^{3}+\frac{1}{288}a^{2}-\frac{7}{480}a \right).\]

We factor out $(a-2)!$, and consider the quadratic 
\[ f_{a}(z) = \frac{a(a-1)}{2}z^{2} + (a-1)\left( \frac{2}{3} a^{3} + \frac{1}{12} a \right) z + \left( \frac{2}{9}a^{6}-\frac{3}{5}a^{5}+\frac{5}{9}a^{4}-\frac{1}{6}a^{3}+\frac{1}{288}a^{2}-\frac{7}{480}a \right).\]

By the quadratic formula (and Maple) we find that the roots of $f_{a}$ are
\[ -\frac{2}{3}a^{2}-\frac{1}{12} \pm \frac{1}{15}\sqrt{170a^{3}-55a^{2}-5a-5}.\]
If $r$ is either of these roots (as we are interested in the limiting behaviour of the roots, either root will do), we find by the Gauss-Lucas theorem that $W_{a}$'s leftmost root $R_{a}$ is to the left of $r$. A straightforward calculation shows that 
\[ \lim _{n->\infty} \frac{r}{n^{2}} = -\frac{1}{24},\]
Thus for any fixed $\varepsilon > 0 $, and sufficiently large $n$, 
\[ R_{a} \leq \left( -\frac{1}{24} + \varepsilon \right) n^{2}.\]
Finally, as the imaginary parts of chromatic roots of $C_{4}(a,a,a,a)$ are the square roots of the roots of $W_{a}$, we find that $C_{4}(a,a,a,a)$, for sufficiently large enough $a$, has a chromatic root with imaginary part at least 
\[ \sqrt{\frac{1}{24} - \varepsilon}.\]

Putting all the pieces together, we have shown:

\begin{theorem}
The growth rate of the maximum imginary part of a chromatic root is linear, that is, there are positive constants $C_{1}$ and $C_{2}$ such that for all sufficiently large $n$,   
\[ C_{1}n \leq \mathrm{maximaginary}(n) \leq C_{2}n. \] \qed
\end{theorem}

In fact the proof shows that any positive constant slightly less than $1/\sqrt{24} \approx 0.2041$ will do for $C_{1}$.

\section{Non-real chromatic roots of almost all graphs}

We now turn to random graphs, and ask, is it more likely that all the chromatic roots are real or not? Our model is the usual Erd\"{o}s-R\'{e}nyi model $G \in {\mathcal G}_{n,p}$, where each edge appears independently  with fixed probability $p$. 
%However, we can say that almost all graphs of order $n$ has $0,1,\ldots, (2-\varepsilon)\log_{2} n$ as chromatic roots (as the clique number of such a graph is greater than $2-\varepsilon)\log_{2} n$). 
Of the $833$ (isomorphism classes of) connected graphs with seven vertices, $273$ of them have chromatic polynomials with only real roots. (For eight vertices the proportion is $1627/11117$.) In this section we prove that for any fixed $p \in (0,1)$, as $n \rightarrow \infty$ almost all random graphs $G \in {\mathcal G}_{n,p}$ have a non-real root.

It is well known that the chromatic polynomial of a graph $G$ of order $n$ and size $m$ is monic, of degree $n$, with integer coefficients of alternating sign. The top coefficients are known (see, for example \cite[p. 31]{dongbook}:
\[ \pi(G,x) = x^{n} - mx^{n-1} + \left( {{m} \choose {2}} - t \right)x^{n-2} - \cdots,\]
where $t$ is the number of triangles (i.e. $K_{3}$'s) in $G$.
The expected number of edges and triangles in a random graph $G \in {\mathcal G}_{n,p}$ are, respectively,
\[ E(K_{2}) = p{{n} \choose {2}} \mbox{ and } E(K_{3}) = p^{3}{{n} \choose {3}}.\]
Chebyshev's inequality for a discrete random variable $X$ states that for any $\lambda > 0 $,
\[ \mbox{Prob}(|X - E(X)| \geq \lambda) \leq \frac{\mbox{Var(X)}}{\lambda^{2}},\]
and it follows that for any $\varepsilon > 0$, that  
\[ \mbox{Prob}(|X - E(X)| \leq \varepsilon |E(X)|) \geq \frac{\mbox{Var(X)}}{\varepsilon^{2}(E(X))^{2}}.\]
Standard techniques can show that for both of of the random variables 
$X = M$, the number of edges, and $T$, the number of triangles, for any graph in ${\mathcal G}_{n,p}$, 
\[\mbox{Var(X)} = o((E(X))^{2}).\]
(For example, writing $T = \sum T_{S}$, where the sum is taken over all subsets of cardinality $3$ of the vertex set $\{1,\ldots,n\}$ and $T_{S}$ is an indicator random variable for whether $S$ induces a triangle, then 
\begin{eqnarray*} 
\mbox{Var}(T) & = & \sum_{S} \mbox{Var}(T_{S}) + \sum_{S^{\prime} \neq S} \mbox{Cov}(T_{S},T_{S^\prime})\\
  & \leq & E(T) + \sum_{S^{\prime} \neq S} \mbox{Cov}(T_{S},T_{S^\prime})\\
  & \leq & E(T) + \sum_{S^{\prime} \neq S, |S| \geq 2, |S^\prime| \geq 2} E(T_{S}T_{S^\prime})\\
  & = & E(T) + 6{{n} \choose {4}}p^{6}\\
  & = & o((E(T))^2),
\end{eqnarray*}
where we have partitioned the pairs of subsets $(S,S^{\prime})$ according to the cardinality of their intersection -- the covariance is $0$ is their intersection is of size $0$ or $1$ and used the fact that $\mbox{Cov}(X,Y) \leq E(XY)$.)
It follows from Chebyshev's inequality that for any fixed $\varepsilon > 0 $, and for almost all graphs $G \in {\mathcal G}_{n,p}$, 

\begin{eqnarray*}
(1-\varepsilon)p{{n} \choose {2}} \leq M \leq (1+\varepsilon)p{{n} \choose {2}},
\end{eqnarray*}
and
\begin{eqnarray*}
(1-\varepsilon)p^{3}{{n} \choose {3}} \leq T \leq (1+\varepsilon)p^{3}{{n} \choose {3}}.
\end{eqnarray*}

Let $G \in {\mathcal G}_{n,p}$. With probability tending to $1$, the values of $M$ and $T$ are
\begin{eqnarray}
M & = & (1 + \varepsilon_{M})p{{n} \choose {2}}, \mbox{ and}\label{quadm}\\
T & = & (1 + \varepsilon_{T})p^{3}{{n} \choose {3}}\label{quadt},
\end{eqnarray}
where $\varepsilon_{M}$ and $\varepsilon_{T},$ are all bounded in absolute value by some fixed but very small $\varepsilon > 0 $, dependent on $p$, that we shall choose shortly.

We again apply the Gauss-Lucas Theorem applied to the $(n-2)$-th derivative of $\pi(G,x)$:
\[ f_{n-2} = (\pi(G,x))^{(n-2)} = (n-2)! \left( \frac{n(n-1)}{2}x^{2} - (n-1)mx + {{m} \choose {2}} - t \right) .\]
By the quadratic formula, $f_{n-2}$ has a non-real root if and only if its discriminant is negative. Substituting in (\ref{quadm}) and (\ref{quadt}), we find that the discriminant of $f_{n-2}/(n-2)!$ is 
\[ p^{2} \left( \frac{1}{3}(\varepsilon_{T}+1)p - \frac{1}{4}(\varepsilon_{M}+1)^2  \right) .\] 
For any $p < 3/4$ we can choose $\varepsilon$ positive but sufficiently close to $0$ to force this discriminant to be negative, and hence for all $p \in (0,3/4)$, $f_{n-2} = (\pi(G,x))^{(n-2)}$ has a nonreal root. The Gauss-Lucas theorem implies the same is true for $\pi(G,x)$.

\begin{figure}[htp]
\begin{center}
\includegraphics[scale=0.6]{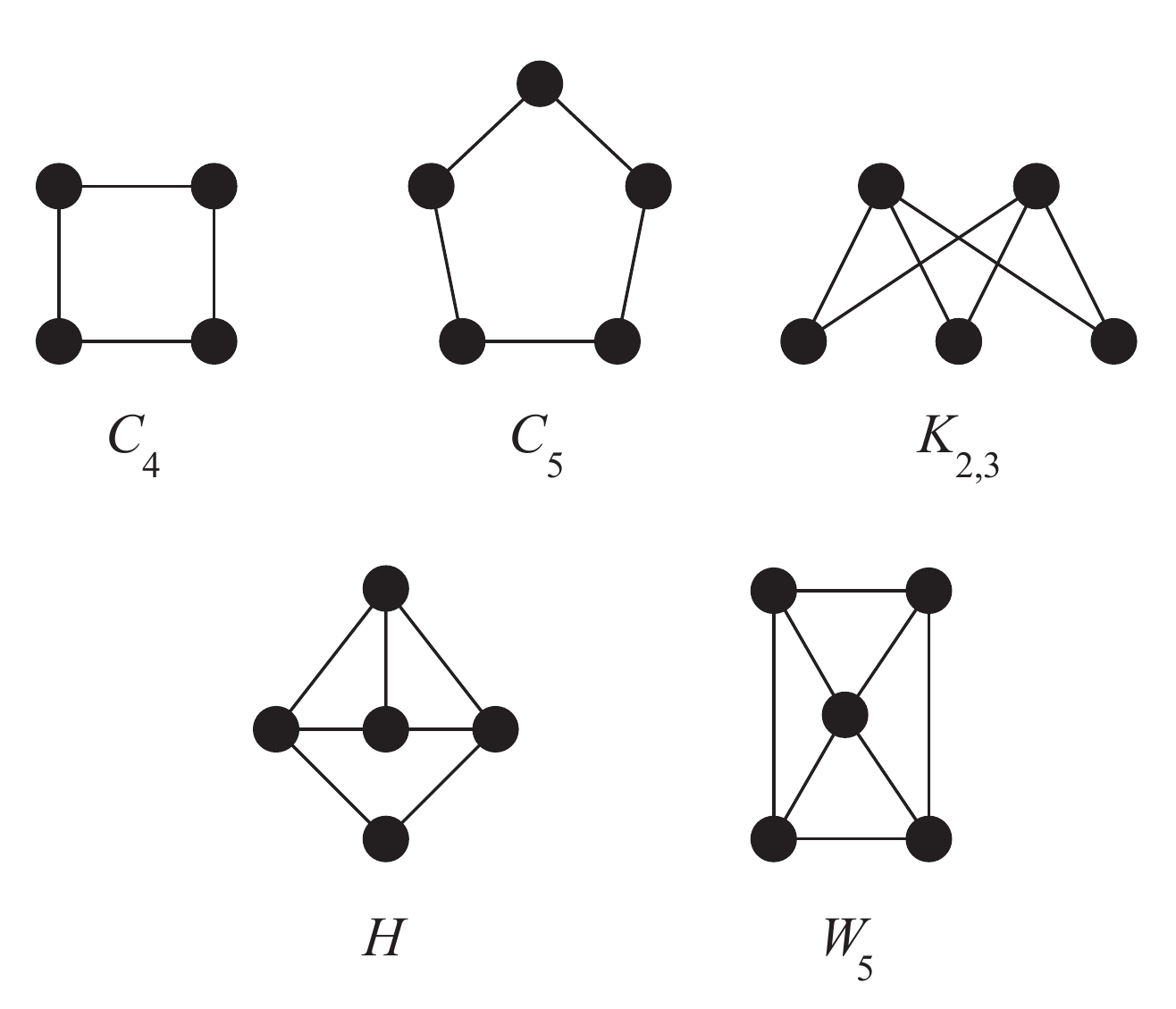}
\caption{Graphs whose counts appear in some of the coefficients in the chromatic polynomial.}
\label{quarticgraphs}
\end{center}
\end{figure}

Now what about $p \geq 3/4$? The argument provided fails, as then in general $f_{n-2}$ has two real roots. We shall need to be more subtle in our argument, and jump from using a quadratic to using a quartic (the use of a cubic provides no assistance here). To do so, we consider the expansion of the chromatic polynomial for the first five terms from the top (again, see \cite[p. 31-32]{dongbook}):
\begin{eqnarray*} 
&&x^{n} - Mx^{n-1} + \left( {{M} \choose {2}} - T \right)x^{n-2} - \left( {{M} \choose {3}} - (M-2)T - IC_{4} + 2nIK_{4} \right) x^{n-2} + \\
&&\biggl( {{M} \choose {4}} - {{M-2} \choose {2}}T + {{T} \choose {2}} -(M-3)IC_{4} -(2M-9)IK_{4} -\\
&&  IC_{5} + IK_{2,3} + 2IH +3IW_{5} -6IK_{5}\biggr) x^{n-3} - \cdots,  
 \end{eqnarray*}
where $IK_{4}$ and $IK_{5}$ are the number of $K_{4}$'s and $K_{5}$'s in $G$, respectively, and $IC_{4}$, $IC_{5}$, $IK_{2,3}$, $IH$ and $IW_{5}$ are the number of {\em induced} $C_{4}$'s, $C_{5}$'s, $K_{2,3}$'s, $H$'s (see Figure~\ref{quarticgraphs}) and $W_{5}$'s (i.e. a wheel of order $5$) in $G$, respectively. 

The expected number of edges, triangles, $K_{4}$-s induced $C_{4}$-s, induced $C_{5}$-s, induced $K_{2,3}$-s, induced $H$-s, induced $W_{5}$-s and $K_{5}$-s in a random graph $G \in {\mathcal G}_{n,p}$ are, respectively:

\begin{equation*}
\begin{aligned}[t]
M & = p{{n} \choose {2}}\\
T & = p^{3}{{n} \choose {3}}\\
IK_{4} & = p^{6}{{n} \choose {4}}\\
IC_{4} & = 3p^{4}(1-p)^{2}{{n} \choose {4}}\\
IC_{5} & = 12p^{5}(1-p)^{5}{{n} \choose {5}} ~~~~~~~~~~~~~~~~~~~~~~~~~~~
\end{aligned}
\begin{aligned}[t]
IK_{2,3} & = 10p^{6}(1-p)^{4}{{n} \choose {5}}\\
IH & = 60p^{7}(1-p)^{3}{{n} \choose {5}}\\
IW_{5} & = 15p^{8}(1-p)^{2}{{n} \choose {5}}\\
K_{5} & = p^{10}{{n} \choose {5}}\\
~ & ~~~~~~~~~~~~~
\end{aligned}
\end{equation*}
%\begin{eqnarray*} 
%M & = & p{{n} \choose {2}},\\
%T & = & p^{3}{{n} \choose {3}},\\
%IK_{4} & = & p^{6}{{n} \choose {4}},\\
%IC_{4} & = & 3p^{4}(1-p)^{2}{{n} \choose {4}},\\
%IC_{5} & = & 12p^{5}(1-p)^{5}{{n} \choose {5}},\\
%IK_{2,3} & = & 10p^{6}(1-p)^{4}{{n} \choose {5}},\\
%IH & = & 60p^{7}(1-p)^{3}{{n} \choose {5}},\\
%IW_{5} & = & 15p^{8}(1-p)^{2}{{n} \choose {5}}, \mbox{ and } \\
%IK_{5} & = & p^{10}{{n} \choose {5}}.\\
%\end{eqnarray*}

Using similar techniques as presented earlier on counting triangles, for all of these random variables and for a graph in ${\mathcal G}_{n,p}$, 
$\mbox{Var(X)} = o((E(X))^{2}),$ and so from Chebyshev's inequality that for any fixed $\varepsilon > 0 $, and for almost all graphs $G \in {\mathcal G}_{n,p}$,
 \begin{eqnarray*}
 (1-\varepsilon)10p^{6}(1-p)^{4}{{n} \choose {5}} & \leq IK_{2,3} \leq & (1+\varepsilon)10p^{6}(1-p)^{4}{{n} \choose {5}},
 ]\end{eqnarray*}
for example. Similar inequalities hold in all other cases.
%\begin{eqnarray*}
%(1-\varepsilon)p{{n} \choose {2}} & \leq M \leq & (1+\varepsilon)p{{n} \choose {2}},
%\\
%(1-\varepsilon)p^{3}{{n} \choose {3}} & \leq T \leq & (1+\varepsilon)p^{3}{{n} \choose {3}},\\
%(1-\varepsilon)p^{6}{{n} \choose {4}} & \leq IK_{4} \leq & (1+\varepsilon)p^{6}{{n} \choose {4}},\\
%(1-\varepsilon)3p^{4}(1-p)^{2}{{n} \choose {4}} & \leq IC_{4} \leq & (1+\varepsilon)3p^{4}(1-p)^{2}{{n} \choose {4}},\\
%(1-\varepsilon)12p^{5}(1-p)^{5}{{n} \choose {5}} & \leq IC_{5} \leq & (1+\varepsilon)12p^{5}(1-p)^{5}{{n} \choose {5}},\\
%(1-\varepsilon)10p^{6}(1-p)^{4}{{n} \choose {5}} & \leq IK_{2,3} \leq & (1+\varepsilon)10p^{6}(1-p)^{4}{{n} \choose {5}},\\
%(1-\varepsilon)60p^{7}(1-p)^{3}{{n} \choose {5}} & \leq IH \leq & (1+\varepsilon)60p^{7}(1-p)^{3}{{n} \choose {5}},\\
%(1-\varepsilon)15p^{8}(1-p)^{2}{{n} \choose {5}} & \leq IW_{5} \leq & (1+\varepsilon)15p^{8}(1-p)^{2}{{n} \choose {5}}, \mbox{ and}\\
%(1-\varepsilon)p^{10}{{n} \choose {5}} & \leq IK_{5} \leq & (1+\varepsilon)p^{10}{{n} \choose {5}}.
%\end{eqnarray*}
 
For $G \in {\mathcal G}_{n,p}$, with probability tending to $1$, the values of the salient graph parameters are
\begin{eqnarray}
M & = & (1 + \varepsilon_{M})p{{n} \choose {2}},\label{quarm}\\
T & = & (1 + \varepsilon_{T})p^{3}{{n} \choose {3}},\label{quart}\\
IK_{4} & = & (1 + \varepsilon_{IK_{4}})12p^{6}{{n} \choose {4}},\label{quark4}\\
IC_{4} & = & (1 + \varepsilon_{IC_{4}})3p^{4}(1-p)^{2}{{n} \choose {4}},\label{quaric4}\\
IC_{5} & = & (1 + \varepsilon_{IC_{5}})12p^{5}(1-p)^{5}{{n} \choose {5}},\label{quaric5}\\
IK_{2,3} & = & (1 + \varepsilon_{IK_{2,3}})10p^{6}(1-p)^{4}{{n} \choose {5}},\label{quarik23}\\
IH & = & (1 + \varepsilon_{IH})60p^{7}(1-p)^{3}{{n} \choose {5}},\label{quarih}\\
IW_{5} & = & (1 + \varepsilon_{IW_{5}})15p^{8}(1-p)^{2}{{n} \choose {5}}, \mbox{ and},\label{quariw5}\\
IK_{4} & = & (1 + \varepsilon_{IK_{4}})p^{10}{{n} \choose {5}},\label{quark5}
\end{eqnarray}
where $\varepsilon_{M},\varepsilon_{T},\varepsilon_{IK_{4}},\varepsilon_{IC_{4}},\varepsilon_{IC_{5}},,\varepsilon_{IK_{2,3}},\varepsilon_{IH},\varepsilon_{IW_{5}}$ and $\varepsilon_{IK_{5}}$ are all bounded in absolute value by some fixed but very small $\varepsilon > 0 $, to be chosen to satisfy some inequalities.

We now apply the Gauss-Lucas Theorem to the $(n-4)$-th derivative of $\pi(G,x)$, which is $(n-4)!$ times

\begin{eqnarray*}
& & \frac{n(n-1)(n-2)(n-3)}{24}x^{4} - \frac{(n-1)(n-2)(n-3)}{6}mx^{3} +
\frac{(n-2)(n-3)}{2}\left( {{m} \choose {2}} - t \right) x^{2} \\
& - & (n-3)\left( {{m} \choose {3}} - (m-2)t - ic_{4} + 2nk_{4} \right) x \\
& + & 
\left( {{m} \choose {4}} - {{m-2} \choose {2}}t + {{t} \choose {2}} -(m-3)ic_{4} -(2m-9)k_{4} - ic_{5} + ik_{2,3} + 2ih +3iw_{5} -6k_{5}\right).
\end{eqnarray*}
 
When a quartic has all real roots is more involved than for a quadratic (or cubic). The {\em discriminant} of a quartic 
\[ g = ax^{4}+bx^{3}+cx^{2}+dx+e\]
is given by
\begin{eqnarray*}
\mbox{Disc}(g) & = & 256a^3e^3-192a^2bde^2-128a^2c^2e^2+144a^2cd^2e-27a^2d^4+144ab^2ce^2-6ab^2d^2e
\\
 & & -80abc^2de+18abcd^3+16ac^4e-4ac^3d^2-27b^4e^2+18b^3cde-4b^3d^3-4b^2c^3e+b^2c^2d^2.
 \end{eqnarray*}
 
If a quartic's discriminant is negative, then the quartic has two distinct real roots and two non-real roots \cite{rees}. 

Substituting in (\ref{quarm})--(\ref{quark5}) to the discriminant of $(n-4)!(\pi(G,x))^{(n-4)}$ (as given in Proposition~\ref{disc4}), we get a polynomial of degree $30$ in $n$, whose leading coefficient we deonte by $lc$. If we set all the various $\varepsilon$'s equal to $0$ in $lc$, we get
\begin{eqnarray*}
& & -(1/93312)p^{21}-(1/186624)p^{20}+(1/124416)p^{19}-(227/80621568)p^{18}-(1/1119744)p^{17}\\
 & & +(5/2985984)p^{16}-(5/2985984)p^{15}+(5/5308416)p^{14}-(1/3538944)p^{13}+(1/28311552)p^{12}.
 \end{eqnarray*}
The largest real root of this polynomial (in $p$) is approximately $0.31564$. As the leading coefficient is negative, it follows that  for $p > 0.32$, this polynomial is negative. The roots of a polynomial depend continuously on its coefficients, so it follows that we can choose $\varepsilon > 0$ so small that if the absolute values of all of $\varepsilon_{M},\varepsilon_{T},\varepsilon_{IK_{4}},\varepsilon_{IC_{4}},\varepsilon_{IC_{5}},\varepsilon_{IK_{2,3}},\varepsilon_{IC_{h}},\varepsilon_{IW_{5}}$ and $\varepsilon_{IK_{5}}$ are at most $\varepsilon$, then $lc$ will be negative, provided that $p > 0.32$ (we ensure that the largest real root is less than $0.32$ and the sign of the leading coefficient of $lc$ remains negative). As $lc$ is the leading coefficient of discriminant of the quartic $(n-4)!(\pi(G,x))^{(n-4)}$), it follows that $(\pi(G,x))^{(n-4)}$ and (by Gauss-Lucas) $\pi(G,x)$ itself has a non-real root. 

By showing that in each of the cases $p < 0.75$ and $p > 0.32$ there is a non-real chromatic root, we have completed our proof of the following.

\begin{theorem}
Let $p \in (0,1)$ be fixed. Then with probability tending to $1$ as $n \rightarrow \infty$, a  graph $G \in {\mathcal G}_{n,p}$ has a non-real chromatic root. \qed
\end{theorem}

%\section{Sigma Polynomials}

%\subsection{New families of $\sigma$-polynomials with nonreal roots}

%Graphs whose complements consist of a $K_{n}$ with a path with $n$ edges hanging off each vertex.

%The complements of circulants $C_{n}(1,2)$ for $n$ odd (they have linear number of triangles).

\section{Concluding Remarks}

We end our discussion with a few questions.

\begin{ques}
Does $\displaystyle{\lim _{n \rightarrow \infty} \frac{\mathrm{maximaginary}(n)}{n}}$ exist? If so, what is its value?
\end{ques}

We have shown that if it does exist, it must be larger than $1/2\sqrt{6} \approx 0.020$. However, even for the ring graphs of order $n$ we considered, the largest imaginary parts seem to approach approximately $0.45n$. And calculations show that the largest imaginary parts of chromatic roots of complete bipartite graphs $K_{n/2,n/2}$ are roughly $0.72n$, which raises an extremal problem. 

\begin{ques}
Which graphs of order $n$ have a chromatic root of largest imaginary part? Is it the complete bipartite graph with (nearly) equal parts?
\end{ques}

We have verified that this is indeed the case for order at most $8$. Finally, while we have shown that almost all graphs ${\mathcal G}_{n,p}$ have a non-real chromatic root, what can be said about the maximum imaginary part?

\begin{ques}
For fixed $p \in (0,1)$, is the maximum imaginary part of a chromatic root of almost all graphs $\Omega(n)$?
\end{ques}

\vskip0.4in
\noindent {\bf \large Acknowledgments:} This research was supported in part by NSERC grants RGPIN 170450-2013 (J.I. Brown) and OGP0105392 (D.G. Wagner).

%\section*{References}

\bibliographystyle{elsarticle-num}
%\bibliography{<your-bib-database>}

\end{document}